\documentclass[letterpaper, 12 pt, onecolumn]{ieeeconf}  

\IEEEoverridecommandlockouts                              
\usepackage{color}

\usepackage{graphics} 
\usepackage{epsfig} 
\usepackage{amsmath} 
\usepackage{amssymb}  

\usepackage[ruled]{algorithm2e}
\usepackage{amsmath}

\renewcommand{\baselinestretch}{1.5}

\newtheorem{theorem}{Theorem}
\newtheorem{remark}{Remark}
\newtheorem{lemma}{Lemma}

\newtheorem{proposition}{Proposition}

\title{\Large \bf
 Minimum Cost Constrained Input-Output and Control Configuration 
Co-Design  Problem: A Structural Systems Approach}

\author{S\'ergio Pequito $^{\dag}$ $\qquad$ Soummya Kar $^{\ddag}$ $\qquad$ George J. Pappas $^{\dag}$
\thanks{This work was supported in part by the TerraSwarm Research Center, one of six centers supported by the STARnet phase of the Focus Center Research Program (FCRP) a Semiconductor Research Corporation program sponsored by MARCO and DARPA, and the NSF ECCS-1306128 grant.}
\thanks{
$^{\ddag}$Department of Electrical and Computer Engineering, Carnegie Mellon University, Pittsburgh, PA 15213}
\thanks{
$^{\dag}$Department of Electrical and Systems Engineering, School of Engineering and Applied Science,
University of Pennsylvania}%
}

\begin{document}

\maketitle
\thispagestyle{empty}
\pagestyle{empty}

\begin{abstract}
In this paper, we study the minimal cost constrained input-output (I/O)  and control configuration \emph{co-design} problem. Given a linear time-invariant \emph{plant}, where a collection of possible inputs and outputs is known a priori, we aim to determine the collection of inputs, outputs and communication  among them incurring in the minimum cost,   such that desired control performance, measured in terms of arbitrary pole-placement capability of the closed-loop system, is ensured. We show that this  problem is NP-hard in general (in the size of the state space). However, the subclass of problems, in which the dynamic matrix is irreducible, is shown to be polynomially solvable and  the corresponding algorithm is presented. In addition, under the same assumption, the same algorithm can be used to solve the minimal cost constrained I/O  selection problem, and the  minimal cost control configuration selection problem, individually.  In order to illustrate the main results of this paper, some simulations are also provided.
\end{abstract}

\section{INTRODUCTION}

Real world systems, such as power systems, public or business organizations, and large manufacturing systems, are often too complex to be tackled by the classical paradigm of mostly centralized decision-making. Such systems are typically characterized by a multitude of decision-makers; due to the distributed nature of the sensing model, in which no decision maker may have a priori access to the entire set of relevant data. Moreover, the communications between the decision-makers may be limited, as is the case in almost all practical networked scenarios. This often rules out the possibility of all-to-all data exchange; hence, centralized data processing and decision-making.  An alternative approach consists of decentralization or decentralized processing, in which the key idea is to equip the individual network decision-makers with autonomous decision-making abilities based on partial system data. 

Clearly, the success of decentralized processing relies on carefully crafting the nature of partial data accessible to the local decision-makers (i.e., the information pattern). Therefore, in this paper, our goal is to identify the critical system locations to be monitored and controlled, and architect the required data exchange between the local components at minimal infrastructure and communication costs, desired closed-loop control performance may be achieved. 
The first set of steps in control systems design thus consist of addressing the following questions \cite{Skogestad04a}:
\begin{itemize}
\item[\textbf{Q1}] Which actuators are required to ensure controllability?
\item[\textbf{Q2}] Which sensors are needed to ensure observability? 
\item[\textbf{Q3}] What is the \emph{information pattern}, i.e.,  which sensors need to supply data to which actuators, such that desired control objectives (for instance, stabilizability) may be ensured?
\end{itemize}

In this paper, we address the \emph{input-output} (I/O) selection  problem (\textbf{Q1}-\textbf{Q2}), and the \emph{control configuration} (CC) selection  problem  (\textbf{Q3}), with the additional constraint that  different actuators, sensors and means-of-communication can incur in different costs. These costs can reflect manufacturing, maintenance and installation costs, or  selection preferences. For instance, consider the selection of phasor measurement units (PMUs)  for state estimation in  power electric grids~\cite{pmustuff}, where sensor parameters such as sampling rate and choice of installation site may dramatically affect the cost, or leader selection problems in which some agents are preferred to others in executing some tasks~\cite{Mesbahi_Egerstedt:2010}. Alternatively, the communication cost may be associated with the use of optic-fiber cable to forward data form sensors to actuators, hence, the cost would depend on its length.

The major focus  of the present paper involves the optimization and qualitative assessment of intrinsic system-theoretic constructs to achieve satisfactory cost-effective decentralized control. Specifically, in this cost-effective decentralized control context, we focus on the \emph{co-design} of sensing-actuation infrastructures and information patterns, i.e., which sensor outputs should communicate or be available to which actuators for feedback.

Utilizing concepts from control theory, graph theory and combinatorial optimization, the major focus of the present paper involves the optimization and qualitative assessment of intrinsic system-theoretic constructs, as well as
 to develop new design/analysis  tools and algorithms that can harness the physical dynamics of such systems to meet the specified large-scale control objectives' guarantees.

\subsection*{Related Work}

Recently, the I/O selection problem have received increasing attention in the literature, especially, since the publication of~\cite{Olshevsky}. In~\cite{Olshevsky}, the \emph{minimal controllability problem} (MCP), i.e., the problem of determining the sparsest input matrix that ensures controllability of a given the system dynamics matrix, was shown to be NP-hard, and  some greedy algorithms provided. Exact solutions to MCP are explored in \cite{MinContNPcomp}, and in \cite{PequitoJournal}, using graph theoretical constructions, the minimal controllability problem is shown to be polynomially solvable for almost all numerical realizations of the dynamic matrix, satisfying a predefined pattern of zeros/nonzeros. Alternatively, in~\cite{Summers,Tzoumas,Clark,Pasqualetti} the configuration of actuators is sought to ensure certain performance criteria; more precisely,  \cite{Summers,Tzoumas,Pasqualetti} focus on optimizing properties of the controllability Grammian,  whereas \cite{Clark} studies leader selection problems, in which leaders are viewed as inputs to the system, and the selection criteria aims to speed up convergence. In addition, in~\cite{Clark,Summers,Tzoumas} the submodularity properties of functions of the controllability Grammian are explored, and design algorithms are proposed that achieve feasible placement with certain guarantees on the optimality gap. The I/O selection problem considered in the present paper differs from the aforementioned problems in the following two aspects: first, the selection of the inputs is restricted to belong to a specific given set of possible inputs, i.e., we study constrained input placement, and, hence, differing from~\cite{Olshevsky,MinContNPcomp,PequitoJournal} in which unconstrained input placement is studied.  Secondly, it contrasts with~\cite{Summers,Tzoumas,Clark,Pasqualetti} in the sense that we do not aim to ensure performance in terms of a function of the controllability Grammian, but we aim to minimize the overall actuation cost, measured in terms of manufacturing/installation/preference costs. Furthermore, instead of optimizing a specific (numerical) system instance, we focus on structural design guarantees that hold for almost all (numerical) system instantiations with a specified dynamic coupling structure. In addition, within the scope of the present problem, we provide optimal solutions under mild assumptions on the dynamics of the system, under very general cost formulations. Note that, if we consider general dynamical systems, even with uniform cost, the problem tackled by us is NP-hard~\cite{NPcompCMIS}. The latter problem is referred to as the minimal constrained input  selection (CMIS) problem, i.e., the problem of determining the minimum number of inputs that ensures structural controllability, which  has been extensively studied, see~\cite{CommaultD13,dionSurvey} and references therein. In a sense, this provides justification for the traditional approaches to solve the I/O selection problem, which include suboptimal methods such as  heuristics, genetic algorithms or relaxations, see for instance \cite{reviewIO,reviewPlacement,  Meyer1994S535, Frecker01042003, WilliamBegg200025}, and references therein. 

Regarding  the CC problem~\cite{reviewIO}, it is worthy to point out that some meaningful advances were recently achieved in terms of determining the numerical gains to achieve closed-loop systems performance, given the existence of  feasible informations patterns, that were accomplished by using convex optimization tools~\cite{Rotkowitz}. More precisely, the  set of feasible solutions is often characterized in terms of a property called \emph{quadratic invariance}, which has been subsequently shown to be necessary and sufficient~\cite{Lessard}, see also~\cite{Mahajan} for a review about the recent developments. Nevertheless, these methods always assume that there exists a feasible information pattern, and no restriction is imposed in terms of the sparsity or the cost incurred by a feasible information pattern. This is one of the goals of the present paper; in a sense, we can use the approach in the present paper to determine feasible information patterns that can, subsequently, be used to determine gains for numerical system instances using the existing tools.  In \cite{Pajic} the design of wireless  control networks is pursued, where given a decentralized plant, modeled as a discrete linear time invariant system equipped with actuators and sensors, the communication topology design between actuators and sensors to achieve decentralized control was posed as a CC selection problem.  Both theoretical and computational perspectives were provided, although the CC selection problem admits a degree of simplification in the discrete time setting. The CC selection problem has also been considered in \cite{Sezer} where a method for determining  the minimum  number of essential inputs and outputs required for  decentralization was provided; however, the characterization does not cope with all cases, see, for instance \cite{Trave} (page 219).   Another work that is related with ours is~\cite{optimumfeedbackpattern}, which studies only the CC selection problem, for a given placement of inputs and outputs, i.e., the inputs and outputs do not incur in any cost. Furthermore, the methodology  in~\cite{optimumfeedbackpattern} is not applicable when some of the  communication costs are taken to be infinite, i.e., when a sensor cannot directly communicate with an actuator, as  is often the case in many large-scale (possibly geographically distributed) systems.

Nevertheless, as also referred in~\cite{reviewIO} there are very few methods that address the I/O and CC problems \emph{jointly}, and to the best of our knowledge none that considers general actuation/sensing and communication costs. Hereafter, we show that the problem is difficult (NP-hard) to solve in general, however, we identify an important subclass that admits polynomial complexity solutions. More precisely, we show that there exist efficient tools to address the I/O and CC selection \emph{co-design} when the dynamics matrix is irreducible. Notice that  this comprises a variety of inter-connected dynamical systems~\cite{interconnectedDyn}, multi-agent networks~\cite{interconnectedDyn,Jadbabaie2003}, or dynamics based in consensus-like protocols~\cite{Dimakis,Jadbabaie2003}, where the irreducibility is essentially ensured by the network \emph{connectivity}.

The closest work to the one presented in here, in the sense that it explores the I/O  and control configuration co-design problem,  is the one in \cite{PequitoJournal}, where the \emph{sparsest} I/O selection and control configuration problems, under the assumption that only the structure of the system dynamics is known but without constraints on the possible inputs/outputs, were addressed. 
 An extension of \cite{PequitoJournal} was presented in~\cite{PequitoCDC13}, where general (possibly heterogeneous) costs to actuate and measure state variables, but with uniform communication or feedback link cost, was considered. However, the I/O cost structure considered in~\cite{PequitoCDC13} is somewhat different from that considered in this paper; in the former, the cost is imposed on the state variables that are to be actuated/observed and not on the specific actuators/sensors as considered in $\mathcal{P}_{1}$ in this paper. In summary, the present work  differs from that presented in~\cite{PequitoJournal,PequitoCDC13} in the following  three major aspects: (i) it considers additional constraints on the possible inputs and outputs used; (ii) the costs depend on the inputs and outputs used to perform a task; and (iii) the communication cost between input-output pairs is arbitrary. 
 \hfill $\circ$
 
 
The main contributions of this paper are as follows: (i)~we show that the minimum cost constrained I/O and control configuration co-design problem is NP-hard; (ii)~we present an efficient algorithm\footnote{The Matlab implementation of the algorithm can be found in https://www.mathworks.com/matlabcentral/fileexchange/49977} (polynomial in the dimension of the state, input and output) to solve it, when the dynamic matrix is irreducible; and (iii)~we show  how  our solution, when the dynamic matrix is irreducible, can be used to solve the minimum cost constrained I/O selection problem, and the minimum cost CC selection problem, individually.


The rest of the  paper is organized as follows. In Section~\ref{probStatement}, we provide the formal problem statement, together with some motivation.   Section~\ref{prelim} reviews some concepts and introduces results  in structural systems theory. Subsequently, in Section~\ref{mainresults} we present the main technical results~(proofs are relegated to the Appendix), followed  by an illustrative example in Section \ref{illustrativeexample}. Conclusions and  discussions on further research are presented in Section \ref{conclusions}.

\section{Problem Statement}\label{probStatement}

In this paper, we  consider a given (possibly large-scale) plant and a collection of  inputs and outputs modeled by
\begin{equation}
\dot x(t) = A x(t) + B u(t),\quad  y(t)=C x(t),\label{LTI}
\end{equation}
where $x(t)=[x_1(t)\ \dots \ x_n(t)]^T\in \mathbb{R}^n$, $u=[u_1(t)\ \dots \ u_p(t)]^T\in \mathbb{R}^p$ and $y=[y_1(t)\ \dots \ y_m(t)]^T \in \mathbb{R}^m$  are the state, input and output, respectively.   In addition, let $\bar A \in \{0,1\}^{n \times n}$ be the binary matrix that represents the structural pattern of $A$, and $\bar B\in \{0,1\}^{n \times p},\bar C\in \{0,1\}^{m \times n}$ the structural patterns of the input and output matrices, respectively. Similarly, let $\bar K\in \{0,1\}^{p\times m}$ be the \emph{information pattern}, where $\bar K_{i,j}=1$ if output $j$ is available to actuator $i$, and zero otherwise.

Further, we aim to ensure that a system achieves the specified large-scale control  \emph{guarantees} when the closed-loop system uses static output feedback, under communication constrains imposed by an information pattern. Specifically, by the careful design of I/O and information pattern infrastructures, we want to ensure that the resulting closed-loop system has no \emph{fixed modes},  as our design guarantee~\cite{Davison}. To this end, denote by $[\bar M]=\{M : M_{ij}=0 \text{ \emph{if} } \bar M_{ij}=0\}$  an equivalence class   of matrices of appropriate dimensions. 
The set of fixed modes of the closed-loop system \eqref{LTI} w.r.t. an information pattern $\bar K$ is given by $\sigma_{\bar K}=\bigcap_{K\in [\bar K]}\sigma(A+BKC)$ (see \cite{Davison}),
where  $\sigma(M)$ denotes the set of eigenvalues of the  matrix $M$. It is known that (see \cite{Davison}) if, for a non-empty symmetric open set $\mathcal{W}\subset \mathbb{C}$, $\sigma_{\bar{K}}\subset\mathcal{W}$ (where $\mathbb{C}$ denotes the set of complex numbers), then there exists a gain $K\in [\bar{K}]$ such that all
the eigenvalues (also known as the poles) of the closed-loop
system $A + BKC$ are in $\mathcal W$.
 Equivalently, we want to ensure that the poles of the closed-loop system can be placed arbitrarily by appropriately tuning the numerical feedback gain parameters under the obtained (designed) information pattern.

Yet, in real-world large-scale systems, more often than not, the exact parameters of the \emph{plant} are not available, or may change over time. Hence, to ensure that the desired closed-loop performance guarantees as discussed above are met, in this paper we adopt a structural design and analysis viewpoint and aim to jointly address the I/O and control configuration (CC) selection such that the closed-loop system has no \emph{structurally fixed modes} (SFMs). The structural version of fixed modes was introduced in \cite{Papadimitriou84asimple}, which, essentially, are the fixed modes attributed to the structural pattern, i.e., location of zeros and nonzeros, of a system, as opposed to fixed modes that originate from a {perfect canceling} of the numerical parameters. Specifically, a structural LTI system $(\bar{A},\bar{B},\bar{C})$ is said to have structurally fixed modes (SFMs) w.r.t. an information pattern $\bar{K}$, i.e., $(\bar A,\bar B,\bar C,\bar K)$ has no SFMs, if for all $A\in [\bar{A}]$, $B\in [\bar{B}]$, $C\in [\bar{C}]$, we have
$\bigcap_{K\in [\bar K]}\sigma(A+BKC)\neq \emptyset$.

Conversely, a structural system $(\bar{A},\bar{B},\bar{C},\bar  K)$ has no structurally fixed modes, if there exists at least one instantiation $A\in [\bar{A}]$, $B\in [\bar{B}]$, $C\in [\bar{C}]$ which has no fixed modes, i.e., $\cap_{K\in [\bar{K}]}\sigma(A+BKC)=\emptyset$. In this latter case, it may be shown (see \cite{sezerFixedModes}) that almost all systems in the sparsity class $(\bar{A},\bar{B},\bar{C})$ have no fixed modes, and, hence, allow pole-placement arbitrarily close to any pre-specified (symmetrical about the real axis) set of eigenvalues. This also justifies our constraint of designing systems with no SFMs in problem~$\mathcal P_1$.

 In summary, we choose the non-existence of SFMs as our design criterion because, informally, it would imply that all LTI systems  represented by $(A,B,C)$ with a given sparsity pattern, i.e., location of zeros/nonzeros  $(\bar{A},\bar{B},\bar{C})$, have no fixed modes, and, hence, would allow pole-placement arbitrarily close to any pre-specified set of eigenvalues.

Thus, in the present paper,  we  address the \emph{minimal cost constrained I/O and control configuration co-design } problem stated as follows.

\subsection*{Problem Statement}

Let $c_u(i)\in\mathbb{R}^n_+$ denote the (non-negative) cost associated with selecting the $i$th input  $u_i$ ($i\in\mathcal I=\{1,\ldots,p\}$), $c_y(j)\in\mathbb{R}^n_+$ denote the (non-negative cost) associated with selecting the $j$th output $y_j$ ($j\in\mathcal J=\{1,\ldots,m\}$), and $c_f((i,j))\in\mathbb{R}^n_+\cup\{\infty\}$, with $(i,j)\in\mathcal I \times\mathcal J$, denote the cost associated with feeding output $j$ to input $i$,  also referred to as communication cost, where $c_f((i',j'))=\infty$ if output $j'$ is not available to input $i'$. This  paper aims to study the following problem.

\vspace{0.2cm}
\noindent $\mathcal P_1$ Find the triple $(\mathcal I^*,\mathcal J^*,\mathcal F^*)$  that solves the following optimization problem:
\begin{align}
\min_{\small\begin{array}{c}\mathcal I\subset \{1,\ldots,p\}\\ \mathcal J\subset \{1,\ldots,m\}\\\mathcal F\subset (\mathcal I\times \mathcal J)\end{array}}&  \sum_{i\in \mathcal I} c_u(i)+\sum_{j\in \mathcal J} c_y(j)+\sum_{(i,j)\in \mathcal F} c_f((i,j))\\
\text{s.t.}\quad & ( \bar A, \bar B(\mathcal I),\bar C(\mathcal J),\bar K(\mathcal F)) \text{ has no SFMs},\notag
\end{align}
where $\bar B(\mathcal I)$ corresponds to the sub-matrix of $\bar B$ comprising the columns with indices in $\mathcal I$, $\bar C(\mathcal J)$ corresponds to the sub-matrix of $\bar C$ comprising the rows with indices in $\mathcal J$, and $\bar K(\mathcal F)_{i,j}=1$ if $(j,i)\in \mathcal F$.

Due to the combinatorial nature of the  I/O and CC co-design problems  (see for instance \cite{reviewIO}) are typically solved using a two-step (generally suboptimal) procedure, namely solving   first the input/output (I/O) selection problem, followed by the control configuration (CC) selection problem. Formally, the structural theory counterparts of these problems are given as follows:

\vspace{0.1cm}
\noindent \textit{Minimum Cost Constrained I/O selection problem}
\vspace{0.1cm}

\noindent $\mathcal P_1^{\text{I/O}}$: Given the structure of the dynamic matrix $\bar A\in\{0,1\}^{n\times n}$, the structure of the input matrix $\bar B\in\{0,1\}^{n\times p}$ and the input costs $c_u(i)$, with $i=1,\ldots,p$,  the \emph{minimum cost constrained input selection}  problem consists in determine $\mathcal I^*$ that solves the following optimization problem:
\begin{align}
\min_{\small\begin{array}{c}\mathcal I\subset \{1,\ldots,p\}\end{array}}& \qquad\qquad \sum_{i\in \mathcal I} c_u(i)\\
\text{s.t.}\qquad & ( \bar A, \bar B(\mathcal I)) \text{ is structurally controllable},\notag
\end{align}
where $\bar B(\mathcal I)$ corresponds to the sub-matrix of $\bar B$ comprising the columns with indices in $\mathcal I$, and a system $(\bar A,\bar B)$ is said to be \emph{structurally controllable} if there exists a controllable  pair   $(A_0,B_0)$ of real matrices, i.e., a system described by these matrices, with zero entries imposed by the zero entries of $(\bar A,\bar B)$. \hfill $\diamond$

\vspace{0.1cm}
\noindent \textit{Minimum Cost Constrained CC selection problem}
\vspace{0.1cm}

\noindent $\mathcal P_1^{\text{CC}}$: Given the structure of the dynamic matrix $\bar A\in\{0,1\}^{n\times n}$, the structure of the input and output matrices $\bar B\in\{0,1\}^{n\times p},\bar C\in\{0,1\}^{m\times n}$, and the communication costs $c_f((i,j))\in\mathbb{R}^n_+\cup\{\infty\}$, with $(i,j)\in(\mathcal I\times \mathcal J) \equiv(\{1,\ldots,p\} \times\{1,\ldots,m\})$,  the \emph{minimum cost constrained control configuration selection}  problem consists in determine $\mathcal F^*$ that solves the following optimization problem:
\begin{align}
\min_{\small\mathcal F\subset (\mathcal I\times \mathcal J)}& \qquad \sum_{(i,j)\in \mathcal F} c_f((i,j))\\
\text{s.t.}\quad & ( \bar A, \bar B,\bar C,\bar K(\mathcal F)) \text{ has no SFMs}.\notag
\end{align}

 \hfill $\diamond$


\section{PRELIMINARIES AND TERMINOLOGY}\label{prelim}

We  start by reviewing some computational complexity concepts~\cite{Garey:1979:CIG:578533}, followed by some concepts related with the study of structural systems theory~\cite{dionSurvey,PequitoJournal}.

\subsection*{Computational Complexity}

A (computational) problem is said to be \emph{reducible in polynomial time} to another if there exists a procedure to transform the former to the latter using a number of operations which is polynomial in the size of its inputs. Such a reduction is useful in determining the  complexity class~\cite{Garey:1979:CIG:578533} a problem belongs to. For instance, recall that a problem $\mathcal{P}$ in  NP  (i.e., the class of non-deterministic polynomial algorithms) is said to be NP-complete if all other NP problems can be polynomially reduced to $\mathcal{P}$~\cite{Garey:1979:CIG:578533}. The set of NP-complete problems is referred to as the NP-complete class.  The following result may be used to prove the NP-completeness of a given problem.

\begin{lemma}[\cite{Garey:1979:CIG:578533}]
    If a problem $\mathcal P_A$ is NP-complete, $\mathcal P_B$ is in NP and $\mathcal P_A$ is reducible in polynomial time to $\mathcal P_B$, then $\mathcal P_B$ is NP-complete.
    \hfill $\diamond$
    \label{NPcompRed}
\end{lemma}

The NP-complete class is used to describe the complexity of decision versions of  problems. For instance,  the following constitutes a decision problem that is particularly relevant in the structural design context: Given $\bar A\in\{0,1\}^{n\times n}$ and $\bar B\in\{0,1\}^{n\times p}$, is there a collection of inputs $\bar B(\mathcal J)$, where $\mathcal J\subset\{1,\ldots,p\}$ and where $\bar B(\mathcal J)$ corresponds to the sub-matrix of $\bar B$ comprising the columns with indices in $\mathcal J$, with $|\mathcal J|=k$ such that $(\bar A,\bar B(\mathcal J))$ is structurally controllable?

Alternatively, it is often natural to consider the optimization versions associated with the decision problems. For instance, the optimization version of the problem stated above aims to determine the minimum $k$ such that the aforementioned property holds. This optimization problem is referred to as the  constrained minimum structural input  selection (CMIS)  problem, given by $\mathcal P_1^{\text{I/O}}$ with uniform non-zero  actuation costs. 
Note that, if a solution to the optimization problem is known, the decision problem is straightforward to solve. Consequently, the optimization problem formulations of NP-complete problems, are  referred to as being NP-hard, since they are at least as difficult as the NP-complete problems; in other words, by solving an instance of the optimization problem (the NP-hard problem), one can obtain a solution to an NP-complete problem.

\subsection*{Graph Theory and Structural Systems}

The following standard terminology and notions from graph theory can be found, for instance, in \cite{PequitoJournal}. Let $\mathcal D(\bar A)=(\mathcal X,\mathcal E_{\mathcal X,\mathcal X})$ be the digraph representation of $\bar{A}$ in (1), where the vertex set $\mathcal X$ represents the set of state variables (also referred to as state vertices) and $\mathcal E_{\mathcal X,\mathcal X}=\{(x_i,x_j): \ A_{ji}\neq 0\}$ denotes the set of edges.  Similarly, we define the following digraphs: $\mathcal D(\bar A,\bar B)=(\mathcal X\cup \mathcal U,\mathcal E_{\mathcal X,\mathcal X}\cup\mathcal E_{\mathcal U,\mathcal X})$ where $\mathcal U$ represents the set of input vertices and $\mathcal E_{\mathcal U,\mathcal X}=\{(u_i,x_j):\ \bar B_{ji}\neq 0\}$; and $\mathcal D(\bar A,\bar B,\bar K,\bar C)=(\mathcal X\cup\mathcal U\cup \mathcal Y,\mathcal E_{\mathcal X,\mathcal X}\cup\mathcal E_{\mathcal X,\mathcal Y}\cup\mathcal E_{\mathcal U,\mathcal X}\cup\mathcal E_{\mathcal Y,\mathcal U})$ denotes the digraph associated with the closed-loop system, where $\mathcal Y$ represents the set of output vertices and $\mathcal E_{\mathcal X,\mathcal Y}=\{(x_i,y_j):\ \bar C_{ji}\neq 0\}$ and the set of feedback edges/links is given by $\mathcal E_{\mathcal Y,\mathcal U}=\{(y_i,u_j):\ \bar K_{ji}\neq 0\}$.

A digraph $\mathcal{D}_s=(\mathcal V_s,\mathcal E_s)$ with $\mathcal V_s\subset \mathcal V$ and $\mathcal E_s\subset \mathcal E$ is called a \textit{subgraph} of $\mathcal{D}$. If $\mathcal V_s=\mathcal V$, $\mathcal{D}_s$ is said to \textit{span} $\mathcal{D}$. Finally, a subgraph with some property $P$ is \emph{maximal} if there is no other subgraph $\mathcal{D}_{s'}=(\mathcal V_{s'},\mathcal E_{s'})$ of $\mathcal{D}$, such that $\mathcal{D}_s$ is a subgraph of $\mathcal{D}_{s'}$, and $\mathcal{D}_{s^{\prime}}$ satisfies property $P$.  A sequence of directed edges $\{(v_1,v_2),(v_2,v_3),\cdots,(v_{k-1},v_k)\}$, in which all the vertices are distinct, is called \textit{an  elementary path} from $v_1$ to $v_k$, as well as a vertex in a digraph with no incoming and outgoing edges (with some abuse of terminology).  A vertex with an edge to itself (i.e., a \emph{self-loop}), or an elementary path from $v_1$ to $v_k$ comprising an additional edge $(v_k,v_1)$, is called a \emph{cycle}. 

In addition, a digraph $\mathcal{D}$ is said to be strongly connected if there exists an elementary path between any  pair of vertices. A \emph{strongly connected component} (SCC) is a maximal subgraph $\mathcal{D}_S=(\mathcal V_S,\mathcal E_S)$ of $\mathcal{D}$ such that for every $v,w \in \mathcal V_S$ there exists a path from $v$ to $w$ and from $w$ to $v$.

Using the above constructions, we can now present some graph theoretical properties that the digraph associated with the closed-loop system $\mathcal D(\bar A,\bar B,\bar C,\bar K)$ must satisfy to avoid the existence of SFMs.

\begin{theorem}[\cite{Pichai3}]
The structural system $(\bar A,\bar B,\bar C)$ associated with \eqref{LTI} has no SFMs w.r.t. an information pattern $\bar K$, if and only if both of the following conditions hold:
\begin{enumerate}
\item[a)] each state vertex $x\in \mathcal{X}$ is contained in a strongly connected component of $\mathcal{D}(\bar A,\bar B,\bar C,\bar K)$ which includes an edge of $\mathcal{E}_{\mathcal{Y},\mathcal{U}}$;
\item[b)] there exists a finite disjoint union of cycles $\mathcal{C}_k=(\mathcal{V}_k,\mathcal{E}_k)$ (subgraph of $\mathcal{D}(\bar A,\bar B,\bar C,\bar K)$) with $k\in\mathbb{N}$ such that $\mathcal{X}\subset \bigcup_{j=1}^{k}  \mathcal V_j $.\hfill $\diamond$
\end{enumerate}
\label{nofixedmodesThm}
\end{theorem}

The conditions in Theorem~\ref{nofixedmodesThm} hold only if the system is both \emph{structurally controllable} and \emph{ structurally observable}\footnote{A system $(\bar A,\bar C)$ is said to be \emph{structurally observable} if there exists an observable  pair   $(A_0,C_0)$ with zero entries enforced by the zero entries in  $(\bar A,\bar C)$.}. More precisely, we have the following result.
 
\begin{proposition}[\cite{Pichai3,PequitoJournal}]
If $(\bar A,\bar B,\bar C,\bar K)$ has no SFMs, then $(\bar A,\bar B)$ and $(\bar A,\bar C)$ are structurally controllable and observable, respectively.
\hfill $\diamond$
\end{proposition}

\subsection{Optimal Assignment Problem}

The optimal assignment problem consists in determining  the collection $\mathcal M^*$ of pairs of indices of a $k\times k$ square matrix $H$ that correspond to the diagonal entries of the matrix $P_1^*HP_2^*$, where $P_1^*$ and $P_2^*$ are permutation matrices such that
\[
(P_1^*,P_2^*)=\arg\min_{P_1,P_2\in\mathcal P} \text{trace}(P_1HP_2),
\]
with $\mathcal P$ denoting the class of all $k\times k$ permutation matrices. 
In what follows, we consider, for the optimal assignment problem,  block matrices $H$ given by
\[
H_{n+p+m}=\left[
\begin{smallmatrix}
H^{x,x}_{n\times n} & H^{x,u}_{n\times p}  &H^{x,y}_{n\times m} \\
H^{u,x}_{p\times n}  & H^{u,u}_{p\times p}  &H^{u,y}_{p\times m} \\
H^{y,x}_{m\times n}  & H^{y,u}_{m\times p}  &H^{y,y}_{m\times m} 
\end{smallmatrix}
\right],
\]
where the labels of the columns associated with the first block are $\{x_1,\ldots,x_n\}$, the second $\{u_1,\ldots, u_p\}$, and the third $\{y_1,\ldots,y_m\}$. Similarly,  the labels of the rows associated with the first block are $\{x_1,\ldots,x_n\}$, the second $\{u_1,\ldots, u_p\}$, and the third $\{y_1,\ldots,y_m\}$. For example, the first row and third column of the $n\times n$ matrix $H^{xx}_{n\times n}$ is indexed by the pair $(x_1,x_3)$,  and the first row and third column of the $n\times p$ matrix $H^{u,x}_{n\times p}$ is indexed by the pair $(u_1,x_3)$. Consequently, any solution to the optimal assignment problem $\mathcal M^*\subset \mathcal L\times \mathcal L$, where $\mathcal L=\{x_1,\ldots,x_n,u_1,\ldots, u_p,y_1,\ldots,y_m\}$, consists of  $n+p+m$ pairs of indices.

Further,  we can associate a matrix $H$ with a digraph representation $\mathcal D(H)=(\mathcal L,\mathcal E)$, where $\mathcal E\subseteq \mathcal L\times \mathcal L$. In particular, a matrix $H$  is \emph{irreducible} if and only if $\mathcal D(H)$ is an SCC. In addition, we  have the following result given by  K{\"o}nig (see Appendix in~\cite{Reinschke:1988}):

\begin{proposition}
Given a square matrix $H$, where $\mathcal D(H)=(\mathcal L,\mathcal E)$, and a solution optimal assignment problem $\mathcal M^*$, then the digraph $\mathcal D=(\mathcal L, \mathcal M^*))$ corresponds to a disjoint collection of cycles that spans $\mathcal D(H)$.
\hfill $\diamond$
\label{KonigProp}
\end{proposition}

\begin{remark}
In the sequel, we will use the solution of  an optimal assignment problem to establish results about 
$\mathcal D(\bar A,\bar B,\bar C,\bar K)$, which is mainly possible due to the labeling we chose to the rows and columns of the  matrix $H$ to be used in the assignment problem, and  consistent with some of the edges in $\mathcal D(\bar A,\bar B,\bar C,\bar K)$. In particular, using Proposition~\ref{KonigProp}, we will be able to ensure that  $\mathcal D(\bar A,\bar B,\bar C,\bar K)$ satisfies condition b) in Theorem~\ref{nofixedmodesThm}.\hfill $\diamond$
\end{remark}

\section{MAIN RESULTS}\label{mainresults}

In this section, we present the main results of the present paper. More precisely, we first show that $\mathcal P_1$ is NP-hard (Theorem~\ref{P1NPcomplete}). However, the subclass of problems, in which the dynamic matrix is irreducible, is shown to be polynomially solvable, using, for instance, Algorithm~1. The correctness and complexity of Algorithm~1 is presented in Theorem~\ref{theoremAlgorithm1}. Finally, we show how Algorithm~1 can be used to solve polynomially  the minimal cost constrained I/O  selection problem  as given in $\mathcal P_1^{I/O}$, and the minimal cost CC problem as given in $\mathcal P_1^{CC}$, if we consider similar assumptions.

We begin by showing that $\mathcal P_1$ is NP-hard. 

\begin{theorem}\label{P1NPcomplete}
The minimal cost constrained I/O  and control configuration co-design problem, given in  $\mathcal P_1$, is NP-hard.
\hfill $\diamond$
\end{theorem}

Nonetheless, the fact that a problem is NP-hard does not preclude the existence of a subclass of problem instances that are easier to solve. In fact, this is the case when we restrict the structure of the dynamics to be irreducible. In Algorithm~1, we present an algorithm to solve $\mathcal P_1$, under the aforementioned constraint. Next, we provide  its proof of correctness and complexity, given in terms of a cubic polynomial in the dimension of the state, input and output.

\begin{algorithm}[!h]

\textbf{Input:} The structural plant matrices $(\bar A,\bar B,\bar C)$ (with $\bar A$ irreducible),  the input costs $c_u(i)$ ($i\in \mathcal I=\{1,\ldots,p\}$), output costs $c_y(j)$ ($j\in \mathcal J=\{1,\ldots,m\}$), and the  communication (feedback) costs  $c_{f}((i,j))$ ($(i,j)\in \mathcal I\times \mathcal J$).

\vspace{0.2cm}
\textbf{Output:} A solution $(\mathcal I^*,\mathcal J^*,\mathcal F^*)$ to $\mathcal P_1$.
\vspace{0.2cm}

\textbf{1.} Let \[
[\mathcal C_{\bar A}]_{i,j}=\left\{\begin{array}{cc}
0,& \text{ if } \bar A_{i,j}=1\\
\infty, & \text{otherwise}
\end{array}
\right. ,\]
and compute the optimal assignment problem $ M'$ considering $\mathcal C_{\bar A}$.

\vspace{0.2cm}

\textbf{2.} \textbf{If} $ M'$ incurs in finite cost, \textbf{then} determine
\[
(i^*,j^*)= \arg \min_{(i,j)\in  \mathcal I\times\mathcal J} c_u(i) +c_y(j) +c_f((i,j)),
\]
and, if $c_u(i^*) +c_y(j^*) +c_f((i^*,j^*))<\infty$, then set $\mathcal I^*=\{i^*\}$, $\mathcal J^*=\{j^*\}$, and $\mathcal F^*=\{(i^*,j^*)\}$, otherwise there is no feasible solution; in particular, since the input and output costs are finite, it follows that there is no feasible information pattern;

\textbf{else} consider the following matrix
\[
\mathcal C^*=\left[\begin{array}{ccc}
\mathcal C_{\bar A}^T & \infty_{n\times p}&  \mathcal C_{\bar C}^T \\
 \mathcal C_{\bar B}^T  & \mathbb{I}^*_{p\times p} & \infty_{p\times m}  \\
\infty_{m\times n} &  \mathcal C_{\bar K}^T& \mathbb{I}^*_{m\times m}
\end{array}\right],
\]
where 
\[ 
[\mathcal C_{\bar B}]_{i,j}=\left\{\begin{array}{cc}
c_u(j),& \text{ if } \bar B_{i,j}=1\\
\infty, & \text{otherwise}
\end{array}
\right. ,
\]
\[ 
[\mathcal C_{\bar C}]_{i,j}=\left\{\begin{array}{cc}
c_y(i),& \text{ if } \bar C_{i,j}=1\\
\infty, & \text{otherwise}
\end{array}
\right.,
\]
and $[\mathcal C_{\bar K}]_{i,j}=c_f((i,j))$, and infinite otherwise.  In addition, 
$\mathbb{I}_{r\times r}^*$ is the $r \times r$ matrix with zero entries in its diagonal and infinity in the off-diagonal entries, and $\infty_{k\times l}$ is the $k\times l$ matrix with its entries set to infinity.

Now,  compute a solution to the optimal assignment problem $M^*$ associated with $\mathcal C^*$. If $M^*$ incurs in finite cost, then set  $\mathcal I^*=\{i \in \mathcal I: \ (u_i,.)\in M^*\}$, $\mathcal J^*=\{j \in \mathcal J: \ (.,y_j)\in M^*\}$, and $\mathcal F^*=\{(i,j) \in \mathcal I\times \mathcal J: \ (y_j,u_i)\in M^*\}$, otherwise there is no feasible solution; in particular, since the input and output  cost are finite, it follows that there is no feasible information pattern;


\caption{Solution to $\mathcal P_1$}	
\label{sol1}
\end{algorithm}

\begin{theorem}\label{theoremAlgorithm1}
Algorithm~\ref{sol1} is correct and its complexity is $\mathcal O((n+m+p)^3)$. \hfill $\diamond$
\end{theorem}

Next, we provide the solution to $\mathcal P_1^{\text{I/O}}$, under the assumption that the structure of the dynamics matrix is irreducible, by resorting to Algorithm~1. Notice that this problem is also NP-hard, since we obtain the CMIS problem by considering uniform non-zero actuation costs. 

\begin{theorem}\label{minCostInput}
If $\bar A$ is irreducible, then $\mathcal P_1^{\text{I/O}}$ is  polynomially  solvable  using  Algorithm 1, when setting  $\bar C=\mathbb{I}_n$ (the $n\times n$ identity matrix),  $c_y(j)=0$, with $j=1,\ldots,n$, and $c_f((i,j))=0$ for $(i,j)\in \mathcal I\times \mathcal J$.

\hfill $\diamond$
\end{theorem}

\begin{remark}
By duality between controllability and observability in LTI systems, Theorem~\ref{minCostInput} readily extends to the minimal cost constrained output  selection, which  consists in determining the minimum number of outputs, given a possible configuration of outputs $\bar C$, that incurs in the minimum cost and ensures structural observability. \hfill $\diamond$
\end{remark}

Similarly, we can solve the minimum cost CC selection problem $\mathcal P_1^{\text{CC}}$, as follows.

\begin{theorem}\label{minCostCC}
If $\bar A$ is irreducible, then $\mathcal P_1^{\text{CC}}$ is polynomially  solvable  using  Algorithm 1, by setting $c_u(i)=0$, with $i\in \mathcal I$, and   $c_y(j)=0$, with $j\in \mathcal J$.
\hfill $\diamond$
\end{theorem}

In the next section, we provide a couple of  examples that illustrate the results attained in this paper.

\section{AN ILLUSTRATIVE EXAMPLE}\label{illustrativeexample}

In this section, we provide two examples where a feasible solution to $\mathcal{P}_{1}$ exists; more precisely,  two different cases in Algorithm~1 are explored.

\subsection*{Example 1}

Let the structure of the dynamics, input and output matrices be given as follows:

\[\bar A=\left[
\begin{smallmatrix}
0& 1&1 &0 & 0& 0  \\
1&0 &0 &1 & 0& 0\\
1 &0&0 &0 & 1& 0\\
0&1 &0&0 & 0& 1\\
0&0 &1&0 & 0& 1 \\
0&0 &0&1 & 1& 0
\end{smallmatrix}\right], \ 
\bar B=\left[
\begin{smallmatrix}
1& 1&0 &0  \\
1&0 &0 &0 \\
0 &1&0 &0 \\
0&0 &1&0 \\
0&0 &0&1  \\
0&0 &1&1 
\end{smallmatrix}\right], \text{ and } \bar C=\left[
\begin{smallmatrix}
1 & 1 & 0 & 0& 0& 0\\
1 & 0 & 1 &0 &0 &0\\
0&0&0&1&1&1
\end{smallmatrix}\right].
\]
In addition, let  the input costs are given by $c_u=[10 \ 10 \ 20 \ 20 ]$, the output costs by $c_y=[ 15 \ 15 \ 50]$, and the communication costs by
\[\bar K=\left[\begin{smallmatrix}
5 & \infty & 25\\
\infty &5 & 25\\
20 &\infty &10\\
\infty &20&10
\end{smallmatrix}\right].
\]
First, notice that $\bar A$ is irreducible, and we can resort to Algorithm~1 to solve $\mathcal P_1$. After we execute Algorithm~1, the solution obtained is $\mathcal I^*=\{1\}$, $\mathcal J^*=\{1\}$ and $\mathcal F^*=\{(1,1)\}$, illustrated in Figure~1-b) by the non-black and non-dashed edges, and incurring in a total cost of $30$.

\subsection*{Example 2}

Let the structure of the dynamics, input and output matrices be given as follows:
\[\bar A=\left[
\begin{smallmatrix}
0& 1&0 &0 & 0  \\
1&0 &1 &1 & 1\\
0 &1&0 &0 & 0\\
0&1 &0&0 & 0\\
0&1 &0&0 & 0 
\end{smallmatrix}\right], \ 
\bar B=\left[
\begin{smallmatrix}
0& 1&0   \\
0&0 &0  \\
0 &0&1  \\
1&0 &0 \\
0&1 &1  
\end{smallmatrix}\right], \text{ and } \bar C=\left[
\begin{smallmatrix}
1 & 0 & 0 & 1& 0\\
0 & 0 & 1 &1 &0 \\
0&0&0&0&1
\end{smallmatrix}\right],
\]
In addition, let  the input costs are given by $c_u=[5 \ 10 \ 10  ]$, the output costs by $c_y=[ 10 \ 10 \ 1]$, and the communication costs by
\[\bar K=\left[\begin{smallmatrix}
10 & 10 & \infty\\
100&\infty & 30\\
\infty&100 &30
\end{smallmatrix}\right].
\]
First, notice that $\bar A$ is irreducible, and we can resort to Algorithm~1 to solve $\mathcal P_1$. The solution obtained is $\mathcal I^*=\{1,2,3\}$, $\mathcal J^*=\{1,2,3\}$ and $\mathcal F^*=\{(2,1),(3,3), (1,2)\}$, illustrated in Figure~2-b) by the non-black and non-dashed edges, and incurring in a total cost of $186$.

\begin{figure}[ht]
\centering
\includegraphics[scale=0.5]{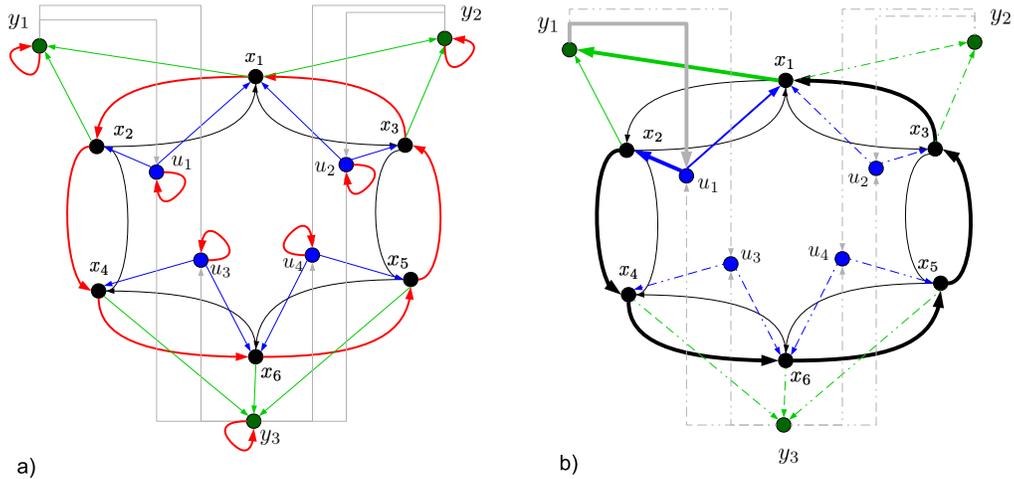}
\caption{In a) we depict the digraph representation associated with $\mathcal C^*$ defined upon the parameters given in Example~1, where a possible solution to the optimal assignment problem is depicted by the red edges. In b) we provide $\mathcal D(\bar A,\bar B(\mathcal I^*),\bar C(\mathcal J^*),\bar K(\mathcal F^*))$ accordingly with the triple $(\mathcal I^*,\mathcal J^*,\mathcal F^*)$ determined by Algorithm~1, where the inputs and outputs with dashed edges have not been selected. In addition, we depict by the bold arrows the disjoint union of cycles in $\mathcal D(\bar A,\bar B,\bar C,\bar K)$ that contains all $x$-vertices, and the SCC containing at least an edge of the form $(y_j,u_i)$, as required by Theorem~1-b) and Theorem~1-a), respectively; hence, Theorem~1 holds, and the closed-loop system $(\bar A,\bar B,\bar C,\bar K)$ has no SFMs. }
\label{example1}
\end{figure}

\begin{figure}[ht]
\centering
\includegraphics[scale=0.5]{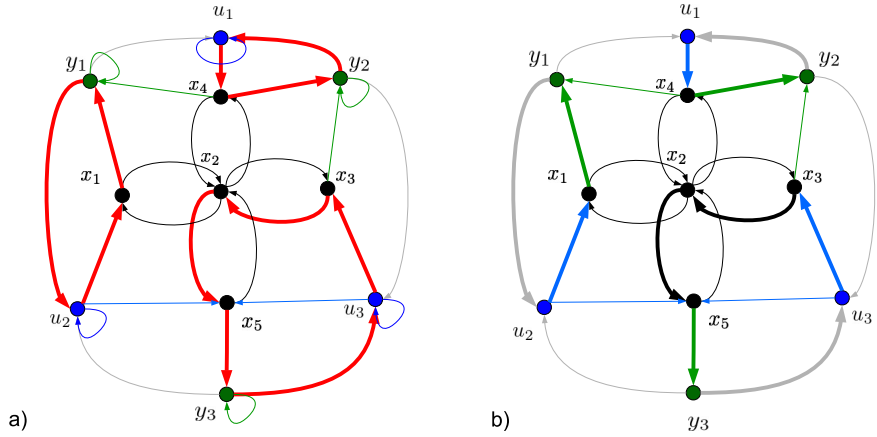}
\caption{In a) we depict the digraph representation associated with $\mathcal C^*$ defined upon the parameters given in Example~2, where a possible solution to the optimal assignment problem is depicted by the red edges. In b) we provide $\mathcal D(\bar A,\bar B(\mathcal I^*),\bar C(\mathcal J^*),\bar K(\mathcal F^*))$ accordingly with the triple $(\mathcal I^*,\mathcal J^*,\mathcal F^*)$ determined by Algorithm~1. In addition, we depict by the bold arrows the disjoint union of cycles in $\mathcal D(\bar A,\bar B,\bar C,\bar K)$ that contains all $x$-vertices, and the SCC containing at least an edge of the form $(y_j,u_i)$, as required by Theorem~1-b) and Theorem~1-a), respectively; hence, Theorem~1 holds, and the closed-loop system $(\bar A,\bar B,\bar C,\bar K)$ has no SFMs. }
\label{example2}
\end{figure}

\section{CONCLUSIONS AND FURTHER RESEARCH}\label{conclusions}

In this paper we have shown that the minimal cost constrained I/O   and control configuration co-design problem  is  an NP-hard problem; hence, efficient algorithms are not likely to exist. Nevertheless, this does not preclude the existence of classes, where it is possible to determine solutions efficiently. In fact, we provided a systematic method with polynomial complexity (in the dimension of the state, inputs and outputs) to jointly solve the input-output and  control configuration selection problem that incurs in a overall minimal cost, under the  assumption that the structure of the dynamics matrix is irreducible.  Future research will consist of determining other subclasses of interest where the current problem can be efficiently solved.

\section*{Appendix}

\subsection*{Proof of Theorem~\ref{P1NPcomplete}}
To prove that $\mathcal P_1$ is NP-hard, we provide a reduction from a known NP-hard problem, the CMIS (see Section~II, problem $\mathcal{P}_1^{\text{I/O}}$), to a particular instance of $\mathcal P_1$ when the costs (input/output/communication) are equal and uniform, which we refer to as $\mathcal P_1'$. Consequently, the corresponding decision problems can be polynomially  reduced to each other, and the result follows by invoking Lemma~\ref{NPcompRed}. First, we notice that it is possible to polynomially verify if  a solution to $\mathcal P_1'$ is feasible, and, consequently, to the corresponding decision problem, see, for instance \cite{optimumfeedbackpattern}; hence, the decision version of $\mathcal P_1'$ is an NP problem. Now, we construct a polynomial reduction from the CMIS problem to $\mathcal P_1'$. Towards this goal, let $\bar A$, $\bar B$ and $c_u(i)$, $i\in \mathcal I=\{1,\ldots,p\}$, in $\mathcal P_1'$ be the same as in the CMIS problem. In addition, let $\bar C=\mathbb{I}_n$ be the $n\times n$ identity matrix,  $c_y(j)=0$ for  $j\in \mathcal J=\{1,\ldots,n\}$ and $c_f((i,j))=0$ for $(i,j)\in \mathcal I\times \mathcal J$. To see that a solution to the proposed problem $\mathcal P_1'$ provides us with a solution to CMIS, recall that a feasible solution to $\mathcal P_1'$, i.e., $(\bar A,\bar B,\bar C,\bar K)$ without SFMs, implies that $(\bar A,\bar B)$ is structurally controllable (see Proposition~1). Now, to see that $B(\mathcal I^*)$ in $\mathcal P_1'$ is also a solution to CMIS, let us assume, by contradiction that it is not. Then, there exists $\mathcal I'$ such that $\bar B(\mathcal I')$ incurs in a lower cost than $\bar B(\mathcal I)$, and such that $(\bar A,\bar B(\mathcal I'))$ is structurally controllable. Now, because $\bar C$ is the identity matrix and $\bar K$ can be full without increasing the cost, it follows that there exists a collection of cycles that comprise the inputs labeled by $\mathcal I'$, as well as a set of labels for the outputs and feedback given by $\mathcal J'$ and $\mathcal F'$, respectively. However, this collection of cycles provides with a solution $\mathcal M'$ to the optimal assignment problem that incurs in lower cost than $\mathcal M^*$, which is a contradiction since we assumed that $\mathcal M^*$ is a solution to the optimal assignment problem.
\hfill $\blacksquare$

\subsection*{Proof of Theorem~\ref{theoremAlgorithm1}}
As discussed in the preliminaries, a solution $\mathcal M^*$ to the optimal assignment problem  of $\mathcal C^*$ always provides a collection of $n+m+p$ pairs of labels of the form $\mathcal M^*\subset \mathcal L\times \mathcal L$, where $\mathcal L=\{x_1,\ldots,x_n,u_1,\ldots, u_p,y_1,\ldots,y_m\}$. Subsequently, by construction of $\mathcal C^*$, if the weight-sum cost of the entries in $\mathcal C^*$ labeled in $\mathcal M^*$ is finite,  then there exists a collection of disjoint cycles  in the digraph representation of $\mathcal C^*$ (see Proposition 2). Now, notice that by construction, only edges associated to finite entries in $\mathcal C^*$ are used; hence, the cycles comprise only edges of the form $(x_i,x_j)$, $(u_i,x_j)$, $(x_i,y_j)$, $(y_j,u_i)$, as well as $(u_i,u_i)$ and $(y_j,y_j)$. Further, the latter edges $(u_i,u_i)$ and $(y_j,y_j)$ do not contribute to ensure either condition a) or b) in Theorem~1, nor are they represented in $\mathcal D(\bar A,\bar B,\bar C,\bar K)$; hence, they can be neglected from the analysis. Subsequently, by noticing  that  if any edge of the form $(u_i,x_j)$, $(x_i,y_j)$, $(y_j,u_i)$ belongs to  $\mathcal M^*$, then so are the other two, otherwise, $\mathcal M^*$ does not comprise a family of cycles with finite weight. Nonetheless, it might be the case that there is no edge of the form $(u_i,x_j)$, $(x_i,y_j)$, $(y_j,u_i)$, that  corresponds to the case where $\mathcal M'$ incurs in finite cost, which implies that there exists a disjoint union of cycles comprising only vertices with labels $\{x_1,\ldots,x_n\}$;  subsequently, all edges of the form $(u_j,u_j)$ and $(y_k,y_k)$ are used in  $\mathcal M^*$. In this case, the weight-sum of $\mathcal M^*$ is equal to zero and, although condition Theorem~1-b) is satisfied, it follows that there is no edge of the form $(y_i,u_j)$ in the SSC containing the state variables in $\mathcal D(\bar A,\bar B,\bar C,\bar K)$. To preclude this case, if $\mathcal M'$ incurs in finite cost, we consider the triple $(\bar B(\{i\}),\bar C(\{j\}),\bar K(\{(i,j)\}))$ that incurs in the smallest  cost; hence, condition Theorem~1-a) is satisfied, and Theorem~1 holds. Alternatively, if $M'$ does not incur in finite cost, then additional edges that are not of the form $(x_i,x_j)$ are required to be associated with a solution to the optimal assignment problem, as result of  Proposition~2. Further, both conditions in Theorem 1 are satisfied: more precisely, a) is satisfied since $\mathcal D(\bar A)$ is an SCC (by assumption), and there must exist at least an edge of the form $(y_j,u_i)$ in $\mathcal M^*$, as well as  in $\mathcal D(\bar A,\bar B(\mathcal I^*),\bar C(\mathcal J^*),\bar K(\mathcal F^*))$  obtained using Algorithm~1. On the other hand, Theorem~1-b) is satisfied by considering,   the edges of a solution to the optimal assignment problem used to define $\mathcal D(\bar A,\bar B(\mathcal I^*),\bar C(\mathcal J^*),\bar K(\mathcal F^*))$.

Finally, we notice that the algorithm's  complexity  is $\mathcal O((n+m+p)^3)$, since it is the complexity of applying the Hungarian algorithm to the optimal assignment problem associated with $\mathcal C^*$; all the remaining steps have lower complexity which renders them moot to the final complexity.
\hfill $\blacksquare$

\subsection*{Proof of Theorem~\ref{minCostInput}}
First, we notice that since $\bar A$ is irreducible, from Theorem~\ref{theoremAlgorithm1}, it follows that Algorithm~1 determines the optimal solution to $\mathcal P_1$. Further,  if $(\bar A,\bar B,\bar C,\bar K)$ has no SFMs, then it is structurally controllable and observable (see Proposition~1).  Because we have that the system is structurally observable and all outputs can be  fed to all inputs without increasing the cost, it follows that the problem consists in determining the  collection of inputs that incurs in the minimum cost; hence, by noticing that this is the same as the $\mathcal P_1^{\text{I/O}}$, the result follows.
\hfill $\blacksquare$





\footnotesize

\bibliographystyle{IEEEtran}
\bibliography{IEEEabrv,acc2015}

\end{document}